\documentclass[journal]{IEEEtran}

\usepackage{amsmath,amssymb,amsthm}
\usepackage{mathrsfs}
\usepackage{bm}
\usepackage{booktabs}
\usepackage{graphicx}
\usepackage{capt-of}
\usepackage{cite}
\usepackage{wrapfig}
\usepackage{xcolor}
\usepackage{tcolorbox}
\usepackage[caption=false,font=footnotesize]{subfig}

\usepackage{algorithm}
\usepackage{algpseudocode}

\usepackage{url}
\usepackage[hidelinks]{hyperref}

\newtheorem{proposition}{Proposition}

\makeatletter
\@ifundefined{c@algorithm}{\newcounter{algorithm}}{}
\makeatother

\begin{document}
\title{A Scalable Bundle Method for Exact Reformulation of SDP in Three-Phase Power Flow Feasibility}

\author{Bohang~Fang,~Lijun~Ding, and~Cong~Chen 
\thanks{\textit{Corresponding authors:~Cong Chen and Lijun Ding.}}%
\thanks{Bohang Fang and Cong Chen (\texttt{\{Bohang.Fang, Cong.Chen\}@dartmouth.edu}) are with Thayer School of Engineering, Dartmouth College, Hanover, NH 03755, USA. Lijun Ding (\texttt{l2ding@ucsd.edu}) is with Department of Mathematics, University of California San Diego, La Jolla, CA 92093.}
}

\maketitle

\begin{abstract}
Power flow feasibility assessment is computationally challenging for unbalanced three-phase distribution networks.
This paper develops a vectorized semidefinite program (SDP) based on the bus injection model (BIM) and reformulates its dual as an exact-penalty problem, enabling us to develop a scalable three-cut proximal bundle method for  feasibility assessment. The proposed bundle method is numerically over $400\times$
faster than   MOSEK with less than
$1/2000$ of its memory; on the decomposed BIM-SDP,  approximately $2\times$ faster
 with  75\% less memory.
\end{abstract}

\begin{IEEEkeywords}
Unbalanced three-phase distribution networks, bundle method, semidefinite program, feasibility.
\end{IEEEkeywords}
\vspace{-0.3cm}
\section{Introduction}
\IEEEPARstart{P}{ower} flow feasibility assessment determines whether a network admits an operating point that satisfies the \emph{power flow equations} and \emph{operational constraints} \cite{cui2019solvability}.
It is essential for hosting-capacity analysis, interconnection studies, and operating-envelope analysis \cite{Molzahn2026DOE}.
While \cite{Molzahn2026DOE} considers single-phase AC networks, distribution networks often require unbalanced three-phase modeling, for which inter-phase coupling and nonconvex power flow equations make feasibility assessment computationally challenging.

To address this challenge, this paper develops a scalable bundle-based solver built on the exact semidefinite program (SDP) relaxation of the bus injection model (BIM) established in \cite{fang2025}, referred to
as the BIM-SDP.
The main contributions beyond \cite{fang2025} are threefold.
First, we derive a vectorized formulation of the BIM-SDP that is suitable for large-scale computation.
Second, we reformulate its dual as an exact-penalty problem and develop a scalable three-cut proximal bundle method, together with an efficient solver for the associated proximal subproblem.
Third, numerical results on large-scale networks demonstrate substantial reductions in runtime and memory relative to CVX/MOSEK baselines.

\vspace{-0.1cm}
\section{Three-Phase Power Network Model}
We follow the BIM for unbalanced three-phase distribution networks \cite{StevenPartI} and derive a {\em vectorized SDP formulation} for AC power flow, supporting our new bundle algorithm in Sec.~\ref{sec:bundle}.

\vspace{-0.3cm}
\subsection{{Primal vectorized BIM-SDP formulation}}\label{subsec:model:compact}

Consider an \(N\)-bus unbalanced three-phase distribution network.
Let \(\mathbf V\in\mathbb C^{3N}\) denote the stacked complex bus-voltage vector, and \(\mathbf Y\in\mathbb C^{3N\times 3N}\) denote the generally sparse network admittance matrix with off-diagonal \(3\times3\) blocks representing line admittances (see \cite[Sec. II]{fang2025} for details). Let \(\mathbf{W}=\mathbf{V}\mathbf{V}^{\mathsf H}\in\mathbb{H}^{3N}\), where \(\mathbb{H}^{3N}\) denotes the set of Hermitian matrices in \(\mathbb{C}^{3N\times 3N}\). 
Under the BIM, with current injection $\mathbf{I}=\mathbf{YV}$, the nodal complex power injection is $\mathbf{P}=\operatorname{diag}(\mathbf{W}\mathbf{Y}^{\mathsf H})\in\mathbb{C}^{3N}$, where $\operatorname{diag}(\cdot)$ extracts the diagonal entries of a matrix as a vector. To achieve a vectorized formulation, we define the linear map \(\mathcal{A}:\mathbb{H}^{3N}\to\mathbb{R}^{18N}\) and the affine map \(\mathbf{m}:\mathbb{R}^{3N}\to\mathbb{R}^{18N}\) as
\[
\begin{aligned}
\mathcal{A}(\mathbf{W})
:=\;&
\Bigl[
\operatorname{Re}(\mathbf{P}),\,
-\operatorname{Re}(\mathbf{P}),\,
\operatorname{Im}(\mathbf{P}),\,
-\operatorname{Im}(\mathbf{P}), \\
&\qquad
\operatorname{diag}(\mathbf{W}),\,
-\operatorname{diag}(\mathbf{W})
\Bigr]^\top, \\ 
\mathbf{m}(\mathbf{u})
:=\;&
\Bigl[
-\mathbf{u}-\overline{\mathbf{p}},\,
\mathbf{u}+\underline{\mathbf{p}},\,
-\overline{\mathbf{q}},\,
\underline{\mathbf{q}},\,
-\overline{\mathbf{v}},\,
\underline{\mathbf{v}}
\Bigr]^\top.
\end{aligned}
\]
where $\rm{Re(\cdot), Im(\cdot)}$ denote the real and imaginary parts and $\overline{\mathbf{p}}, \underline{\mathbf{p}}, \overline{\mathbf{q}}, \underline{\mathbf{q}}, \overline{\mathbf{v}}, \underline{\mathbf{v}}$ are the active power, reactive power, and {\em squared-voltage} limits.
Let \(\mathbf u := (u_j^\phi)_{j=1,\ldots,N,\;\phi\in\{a,b,c\}}
\in \mathbb R^{3N}\) denote the stacked active power injection vector,
where \(u_j^\phi\) is the injection at bus \(j\) on phase \(\phi\).

For a given $\mathbf{u}$, we formulate the following vectorized BIM-SDP for feasibility assessment:
\vspace{-0.15em}
\begin{subequations}\label{eq:FP_compact_slack_SDP}
\begin{align}
\min_{\mathbf{W}\succeq 0,\;\mathbf{z}\geq 0}
\;& \beta\cdot\mathbf{1}^\top\mathbf{z} + \operatorname{tr}(\mathbf{C}\mathbf{W})
&&
\label{eq:FP_compact_slack_SDPa}\\
\text{s.t.}\quad
& \mathcal{A}(\mathbf{W})+\mathbf{m}(\mathbf{u}) \leq \mathbf{z},
&& :\ \mathbf{y}
\label{eq:FP_compact_slack_SDPb}\\
& [\mathbf{W}]_{11}=\mathbf{V}_{1}\mathbf{V}_{1}^{\mathsf H},
&& :\ \bm{\Gamma}
\label{eq:FP_compact_slack_SDPc}
\end{align}
\end{subequations}
Here, $\mathbf{z}$ is a nonnegative slack variable that captures constraint violation, and $\operatorname{tr}(\mathbf{C}\mathbf{W})$, with $\mathbf{C}:=\frac{1}{2}(\mathbf{Y}+\mathbf{Y}^{\mathsf H})$, penalizes the network power loss.
Constraint \eqref{eq:FP_compact_slack_SDPc} fixes the slack-bus voltage, where $\mathbf V_1=[1,e^{-\mathrm{i}\frac{2\pi}{3}},e^{\mathrm{i}\frac{2\pi}{3}}]^\top$ is the prescribed slack-bus voltage phasor in per unit.
The power loss penalty is included in \eqref{eq:FP_compact_slack_SDPa} because, under radial topology and sufficient small $\beta$, the SDP relaxation admits a rank-one optimal solution, i.e., $\operatorname{rank}(\mathbf{W}^\star)=1$, as proved in \cite[Theorem 2]{fang2025}.
Therefore, if $\mathbf 1^\top \mathbf z^\star=0$, a bus-voltage vector $\mathbf{V}$ can be recovered from $\mathbf W^\star$ by rank-one factorization \cite{StevenPartI}, and it satisfies all network constraints. 
Hence, $\mathbf u$ \emph{passes} the feasibility assessment.

Let 
$\mathbf{y}\in\mathbb{R}^{18N}$ and $\bm{\Gamma}\in\mathbb{H}^3$ denote the dual variables associated with
\eqref{eq:FP_compact_slack_SDPb} and \eqref{eq:FP_compact_slack_SDPc}, respectively.
We introduce this vectorized SDP formulation to support the dual bundle method in Sec.~\ref{sec:bundle}. 
In particular, it enables efficient \emph{sparse} matrix--vector products in the leading-eigenvalue/vector computation, which is a key component of our bundle algorithm.

\vspace{-0.2cm}
\subsection{Dual SDP and Exact Penalty Reformulation}

We next derive a vectorized dual SDP formulation of \eqref{eq:FP_compact_slack_SDP} and reformulate it into an exact-penalty problem, which forms the basis for the proposed three-cut proximal bundle method. 
The step-by-step derivation and computational advantages of the dual SDP
are presented in Online Appendix~\ref{appen:dual:FPbeta}.

Let $\mathcal{B}:\mathbb{H}^{3N}\to\mathbb{H}^3$ denote the block-extraction
operator $\mathcal{B}(\mathbf{W}) := [\mathbf{W}]_{11}$, and define
$\mathbf{M}_{1} := \mathbf{V}_{1}\mathbf{V}_{1}^{\mathsf H}$.
Let $\mathcal{A}^*$ and $\mathcal{B}^*$ denote the adjoint operators of
$\mathcal{A}$ and $\mathcal{B}$. In particular,
$\mathcal{B}^*(\bm{\Gamma})$ places $\bm{\Gamma}$ in the top-left
$3\times 3$ block of a $3N\times 3N$ Hermitian matrix, with zeros elsewhere.
Then, the dual of \eqref{eq:FP_compact_slack_SDP} is:
\begin{subequations}\label{eq:dual_FP_beta}
\begin{align}
\min_{\mathbf{y},\bm{\Gamma}}\quad
& -\,\mathbf{m}(\mathbf{u})^\top \mathbf{y}
+ \operatorname{tr}\!\left(\bm{\Gamma}\mathbf{M}_{1}\right)
\label{eq:dual_FP_betaa}\\
\text{s.t.}\quad
& \mathbf{0}\le \mathbf{y} \le \beta\cdot\mathbf{1},\;\text{where $\mathbf{1}$ is the all one vector,}
\label{eq:dual_FP_betab}\\
& \mathbf{H}(\mathbf{y},\bm{\Gamma})
:= \mathbf{C} + \mathcal{A}^*(\mathbf{y}) + \mathcal{B}^*(\bm{\Gamma}) \succeq 0,
\label{eq:dual_FP_betac}
\end{align}
\end{subequations}
The positive semidefinite (PSD) constraint \eqref{eq:dual_FP_betac} makes \eqref{eq:dual_FP_beta} difficult to solve.  
We reformulate \eqref{eq:dual_FP_beta} as a convex exact penalty formulation in Proposition \ref{prop:dual_exact_penalty}
by defining 
\begin{equation}\label{eq:set}
\mathcal{X}:=\{(\mathbf{y},\bm{\Gamma})\mid \mathbf{0}\le \mathbf{y}\le \beta\cdot\mathbf{1},\ \bm{\Gamma}\in\mathbb{H}^3\}, 
\end{equation}
\[
\phi(\mathbf{y},\bm{\Gamma}):=\big[\lambda_{\max}(-\mathbf{H}(\mathbf{y},\bm{\Gamma}))\big]_+,
\ \text{where} \
[a]_+:=\max\{a,0\}.
\]
Here,  \(\lambda_{\max}(\cdot)\) denotes the largest eigenvalue.
\begin{proposition}[Exact penalty reformulation]
\label{prop:dual_exact_penalty}
Suppose 
\eqref{eq:FP_compact_slack_SDP} has a unique  optimal matrix solution $\mathbf{W}^\star$. Then, for any
$\alpha>\operatorname{tr}(\mathbf{W}^\star)$, the optimal solution set of \eqref{eq:dual_FP_beta} equals that of 
\begin{equation}\label{eq:exact_penalty_obj}
\min_{(\mathbf{y},\bm{\Gamma})\in\mathcal{X}}
f(\mathbf{y},\bm{\Gamma})
:= -\mathbf{m}(\mathbf{u})^\top\mathbf{y}
+ \operatorname{tr}(\bm{\Gamma}\mathbf{M}_{1})
+ \alpha\,\phi(\mathbf{y},\bm{\Gamma}).
\end{equation}
\end{proposition}
The proof, following \cite{ding2023revisiting}, is provided in the online Appendix \ref{appen:proofproposition1}. 
Constraint \eqref{eq:dual_FP_betac} is absorbed into the objective as a nonsmooth exact penalty, yielding a convex minimization over the simple set $\mathcal{X}$.  After solving \eqref{eq:exact_penalty_obj}, we could recover $\mathbf{H}^\star(\mathbf{y}^\star, \bm{\Gamma}^\star)$ via \eqref{eq:dual_FP_betac}. By the KKT complementarity $\mathbf{H^\star}\mathbf{W}^\star=\mathbf{0}$, the primal solution $\mathbf{W}^\star$ can be obtained for the feasibility of $\mathbf{u}$.

\vspace{-0.2cm}
\section{Three-Cut Proximal Bundle Method}\label{sec:bundle}

To solve \eqref{eq:exact_penalty_obj}, we exploit the particular nonsmooth structure of 
$f$ in \eqref{eq:exact_penalty_obj},  
and develop a new \emph{three-cut proximal bundle method} for solving \eqref{eq:exact_penalty_obj}. 
In the following, with $\mathcal{X}$ defined in \eqref{eq:set}, we denote a point in $\mathcal X$ by
$\mathbf x:=(\mathbf y,\bm\Gamma)\in\mathbb{R}^{18N}\times \mathbb{H}^3$,
and define the inner product
$\langle \mathbf x_1,\mathbf x_2\rangle
:= \mathbf y_1^\top \mathbf y_2+\operatorname{tr}(\bm\Gamma_1\bm\Gamma_2)$. 

As an overview, our method is iterative and produces two  sequences of points in $\mathcal{X}$ denoted as $\mathbf{x}^k$ and $\mathbf{z}^k$, referred to as \emph{bundle center} and \emph{trial points} later.  We maintain a set (a.k.a. \emph{bundle}) of three functions, called \emph{cuts}, that globally lower bounds $f$ in \eqref{eq:exact_penalty_obj}. At each iteration $k$, (i) we construct an approximation of $f$ (Sec.~\ref{sec:Three-Cut}) based on  the three-cut bundle model; (ii) we solve a proximal subproblem based on the model with center $\mathbf{x}^k$ (Sec.~\ref{subsec:proximal:subproblem}) and obtain its solution $\mathbf{z}^{k+1}$; (iii) based on $\mathbf{z}^{k+1}$, we update the bundle, set the next center $\mathbf{x}^{k+1}$, and repeat until the stopping criterion is met(Sec.~\ref{subsec:bundle:procedure}).

\vspace{-0.3cm}
\subsection{Approximation with Three-Cut Bundle Model}\label{sec:Three-Cut}

\begin{wrapfigure}{l}{0.45\columnwidth}
    \centering
    \includegraphics[width=\linewidth,
    trim=0 0.75cm 0.3cm 1.3cm, clip]{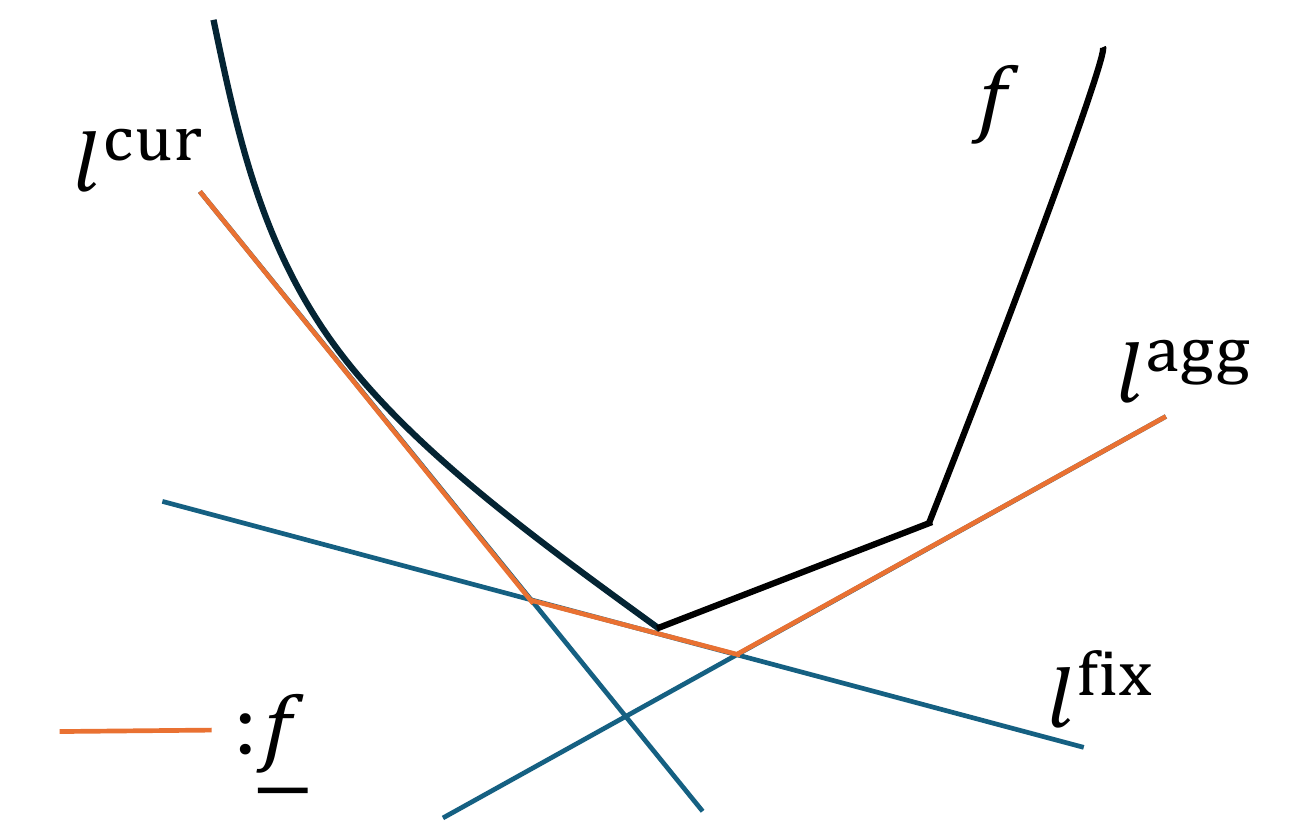}
    \vspace{-0.75cm}
    \caption{Three-cut bundle \eqref{eq:three_cut_model} is a piecewise-affine lower approximation to nonsmooth $f$ in \eqref{eq:exact_penalty_obj}.}
    \label{fig:cuttingplanes}
    \vspace{-0.5cm}
\end{wrapfigure}

To start, suppose for the moment that for each iteration $k\geq1$, we have the three  \emph{cuts}:  three \emph{affine} functions $\ell_k^{\rm fix}$, $\ell_k^{\rm cur}$, and $\ell_k^{\rm agg}$, from $\mathbb{R}^{18N}\times \mathbb{H}^3$ to $\mathbb{R}$, such that they are all \emph{global lower bounds} $f$.
Then, based on the \emph{bundle} $\{\ell_k^{\rm fix},\;\ell_k^{\rm cur},\;\ell_k^{\rm agg}\}$, the {\em three-cut bundle model} $\underline f_k(\mathbf{x})$ defined in \eqref{eq:three_cut_model} is convex and a lower approximation of \eqref{eq:exact_penalty_obj}, see Fig.~\ref{fig:cuttingplanes} for an illustration.
\begin{equation}\label{eq:three_cut_model}
\begin{array}{cc}
\underline f_k(\mathbf{x})
:= &
\!\max\bigl\{
\ell_k^{\rm fix}(\mathbf{x}),\ell_k^{\rm cur}(\mathbf{x}),
\ell_k^{\rm agg}(\mathbf{x})
\bigr\}.
\end{array}
\end{equation}
\noindent 
The rationale for using three cuts is to capture the nonsmoothness of $f$
while ensuring convergence: the fixed cut $\ell_k^{\rm fix}
= -\mathbf{m}(\mathbf{u})^\top\mathbf{y}
+ \operatorname{tr}(\bm{\Gamma}\mathbf{M}_{1})$ captures the linear part of $f$ and 
the current cut $\ell_k^{\rm cur}$ captures the complicated eigenvalue part
near the bundle center, and the maximum of them captures the nonsmoothness of $f$ due to the max and eigenvalue in $\phi$,
while $\ell_k^{\rm agg}$ serves as a global stabilizer
of the algorithm. We leave the detailed initialization and updates of the cuts
to Sec.~\ref{subsec:bundle:procedure}.

\vspace{-0.3cm}
\subsection{Trial Point Update via a Proximal Subproblem}
\label{subsec:proximal:subproblem}

Given the $k$-th bundle center $\mathbf x^k$ and the model
$\underline f_k$, we solve the following proximal subproblem and take
its optimal $\mathbf x$-solution as the next trial point
$\mathbf z^{k+1}$:
\begin{equation}\label{eq:prox_subproblem_epi}
\begin{aligned}
\min_{\mathbf x\in\mathcal X,\,r}\quad
& r+\frac{\rho}{2}\|\mathbf x-\mathbf x^k\|^2 \\
\mathrm{s.t.}\quad
& \ell_k^i(\mathbf x)=a_k^i+\langle \mathbf h_k^i,\mathbf x\rangle \le r,
\quad i=1,2,3,
\end{aligned}
\end{equation}
where
$(\ell_k^1,\ell_k^2,\ell_k^3)
:=
(\ell_k^{\rm fix},\ell_k^{\rm cur},\ell_k^{\rm agg})$,
$a_k^i\in\mathbb R$,
$\mathbf h_k^i\in\mathbb R^{18N}\times\mathbb H^3$, and
$r$ is the epigraph variable for
$\underline f_k(\mathbf x)=\max_i\ell_k^i(\mathbf x)$.
The parameter $\rho>0$ controls the proximal regularization.

A \emph{key innovation} of this paper is that we develop Algorithm~\ref{alg:prox_solver} to solve the \emph{dual} of \eqref{eq:prox_subproblem_epi} \emph{efficiently}, and quickly recover the solution to \eqref{eq:prox_subproblem_epi}.
We prove Proposition~\ref{prop:prox_dual} and the procedure in Algorithm~\ref{alg:prox_solver}, with  details in Appendices \ref{app:proof_prop_prox_dual}-\ref{app:prox_cases}.

\begin{proposition}\label{prop:prox_dual}
Algorithm \ref{alg:prox_solver} produces a dual optimal solution $\bm \theta^k\in\Delta_3$
of \eqref{eq:prox_subproblem_epi}, 
where $\Delta_3:=\{\bm\theta\in\mathbb R^3:\bm\theta \geq 0, \;\mathbf 1^\top\bm\theta=1\}$, and the optimizer of
\eqref{eq:prox_subproblem_epi} is uniquely
\begin{equation}\label{eq:trial_point_general_short}
\mathbf z^{k+1}
=
\Pi_{\mathcal X}\!\left(
\mathbf x^k-\frac{1}{\rho}\sum_{i=1}^3 \theta_i^k \mathbf h_k^i
\right),
\end{equation}
where $\Pi_{\mathcal X}(\cdot)$ is the Euclidean projection\footnote{Based on $\cal X$ defined in \eqref{eq:set}, $\Pi_{\mathcal X}$ acts only
on the $\mathbf y$-block
by clipping each constrained entry to $[0,\beta]$, while leaving the
$\bm\Gamma$-block unchanged.
} onto $\mathcal X$. 
\end{proposition}

While the dual of \eqref{eq:prox_subproblem_epi} is a quadratic program (detailed in \eqref{eq:prox_dual}), using a generic quadratic-program solver would easily incur high computation cost. 
Instead, our Algorithm~\ref{alg:prox_solver} exploits the geometry of $\Delta_3$ and solves the dual problem 
{\em sequentially} over {\em three cases}: vertices, edges, and the interior of $\Delta_3$, in near $\mathcal{O}(N)$ time. Moreover, its sequential design stops at the first successful case, thereby reducing computation in practice.

\begin{figure}[t]
\centering
\refstepcounter{algorithm}
\label{alg:prox_solver}

\begin{minipage}{0.99\columnwidth}
\hrule height 0.8pt
\vspace{0.25em}

\noindent\textbf{Algorithm \thealgorithm:} Proximal Subproblem Solver

\vspace{0.15em}
\hrule height 0.4pt
\vspace{0.25em}

\noindent
{\footnotesize
\textit{Note:} The solver proceeds \emph{sequentially} with early termination:
it terminates once Case~1 returns a solution. Case~3 is entered only
if the first two cases did not return a solution.
\(\mathbf e_i\in \mathbb R^3\) is the \(i\)-th basis vector.
}

\vspace{0.15em}
\noindent
\textbf{Input:}
\(\mathbf x^k\), \(\rho\), and
\(\{\ell_k^i(\mathbf x)=a_k^i+\langle\mathbf h_k^i,\mathbf x\rangle\}_{i=1}^3\).

\vspace{0.15em}
\noindent
\textbf{Case 1 test (vertex):}
For each \(i\in\{1,2,3\}\), compute
\(\mathbf z_i=
\Pi_{\mathcal X}\!\left(\mathbf x^k-\rho^{-1}\mathbf h_k^i\right)\).
If
\(\ell_k^i(\mathbf z_i)=\max_{j=1,2,3}\{\ell_k^j(\mathbf z_i)\}\)
\textbf{return}
\((\mathbf z^{k+1},\bm\theta^k)=(\mathbf z_i,\mathbf e_i)\).

\vspace{0.15em}
\noindent
\textbf{Case 2 test (edge):}
For each edge \((i,j)\) of \(\Delta_3\), parameterize
\(\bm\theta(s)=s\mathbf e_i+(1-s)\mathbf e_j\), \(s\in[0,1]\),
and solve the associated one-dimensional problem by the sweep procedure in Appendix~\ref{app:prox_cases:case2} to obtain
\(s^\star\in\mathbb R\). Let
\(\mathbf z_i=
\Pi_{\mathcal X}\!\left(
\mathbf x^k-\frac{1}{\rho}\sum_{q=1}^3
\theta_q(s^\star)\mathbf h_k^q
\right)\).
If
\(\max\{\ell_k^i(\mathbf z_i),\ell_k^j(\mathbf z_i)\}
\geq \ell_k^l(\mathbf z_i)\)
where \(l\notin\{i,j\}\),
\textbf{return}
\((\mathbf z^{k+1},\bm\theta^k)=(\mathbf z_i,\bm\theta(s^\star))\).

\vspace{0.15em}
\noindent
\textbf{Case 3 solve (interior):}
If neither Case~1 nor Case~2 succeeds, solve the interior of \(\Delta_3\) by the semismooth Newton method in Appendix~\ref{app:prox_cases:case3} to obtain
\(\bm{\theta}^k\), and recover \(\mathbf z^{k+1}\) from
\eqref{eq:trial_point_general_short}.

\vspace{0.15em}
\noindent
\textbf{Output:}
\(\mathbf z^{k+1}\) and \(\bm\theta^k\).

\vspace{0.25em}
\hrule height 0.8pt
\end{minipage}
\vspace{-1.0em}
\end{figure}

\vspace{-0.2cm}
\subsection{Initialization, Update, and Step Acceptance}\label{subsec:bundle:procedure}

Denote $\mathbf{x}^0$ as the initial point, and set $\ell_0^{\rm fix}(\mathbf{x})=\ell_0^{\rm cur}(\mathbf{x})=
\ell_0^{\rm agg}(\mathbf{x}) =  -\mathbf{m}(\mathbf{u})^\top\mathbf{y}
+ \operatorname{tr}(\bm{\Gamma}\mathbf{M}_{1})$. Once the trial point $\mathbf z^{k+1}$ and the multiplier vector $\bm{\theta}^k$ are obtained, we define the weighted slope 
\vspace{-0.3cm}
\begin{equation}\label{eq:agg_subgrad}
\bar{\mathbf h}_k := \sum_{i=1}^3 \theta_i^k \mathbf h_k^i
\in \partial \underline f_k(\mathbf z^{k+1}),
\end{equation}
where \(\partial\) denotes the convex subdifferential.
Define $f^{\lambda}(\mathbf{y},\bm{\Gamma})
:= -\mathbf{m}(\mathbf{u})^\top\mathbf{y}
+ \operatorname{tr}(\bm{\Gamma}\mathbf{M}_{1})
+ \alpha\lambda_{\max}(-\mathbf{H}(\mathbf{y},\bm{\Gamma}))$. We then 
update the {\em three-cut bundle} by
\begin{equation}\label{eq:agg_cut_update_letter}
\begin{aligned}
\ell_{k+1}^{\mathrm{agg}}(\mathbf x)
&=\underline{f}_k(\mathbf z^{k+1}) + \left\langle
\bar{\mathbf h}_k,\, \mathbf x-\mathbf z^{k+1}
\right\rangle,\\
{\ell_{k+1}^{\mathrm{cur}}(\mathbf x)}
&= f^{\lambda}(\mathbf z^{k+1}) + \langle
\mathbf{g}^{k+1},\mathbf{x}-\mathbf{z}^{k+1}\rangle,
\end{aligned}
\end{equation}
and $\ell_{k+1}^{\mathrm{fix}}=\ell_{k}^{\mathrm{fix}}$ is fixed. 
Here, $\mathbf g^{k+1}\in \partial f^{\lambda}(\mathbf z^{k+1})$.
Specifically, let $\mathbf{v}^{k+1}$ be a unit eigenvector of $-\mathbf{H}(\mathbf{z}^{k+1})$ corresponding to $\lambda_{\max}(-\mathbf{H}(\mathbf{z}^{k+1}))$ and $\mathbf{Q}^{k+1}=\mathbf{v}^{k+1}(\mathbf{v}^{k+1})^{\mathsf H}$, then
\begin{equation*}
\mathbf{g}^{k+1}
=
\Big(
-\,\mathbf{m}(\mathbf{u})-\alpha\,\mathcal{A}(\mathbf{Q}^{k+1}),\;
\mathbf{M}_{1}-\alpha\,\mathcal{B}(\mathbf{Q}^{k+1})
\Big). 
\end{equation*}
Note that this step is efficient due to {\em the vectorized formulation} in \eqref{eq:FP_compact_slack_SDP} and {\em sparsity} of $\mathbf{Y}$.\footnote{Indeed, the sparsity of $\mathbf{Y}$ enables (i) quick application of $\mathcal{A}$ and $\mathcal{B}$ to rank-one matrices and (ii) the fast matrix-vector operations associated with $\mathcal A^*$ and $\mathcal B^*$ in computing the eigenvector $\mathbf{v}^{k+1}$ using the Lanczos method. } 
To assess whether the trial point yields sufficient improvement, we define the predicted decrease $\Delta_k
:=
f(\mathbf x^k)-\underline f_k(\mathbf z^{k+1})\geq0$, since \(\underline f_k\le f\).
Given $\eta\in(0,1)$, we accept a descent step if
\vspace{-0.1cm}
\begin{equation}\label{eq:descent_step_test_letter}
f(\mathbf z^{k+1})
\le
f(\mathbf x^k)-\eta\,\Delta_k.
\end{equation}
If \eqref{eq:descent_step_test_letter} holds, the trial point is accepted as the next center; otherwise, the current center is retained, corresponding to a null step. 
In implementation, the algorithm is terminated when $\Delta_k \le \varepsilon$, where $\varepsilon>0$ is a prescribed tolerance.
The overall update is summarized in Algorithm~\ref{alg:bundle}. 
\vspace{-0.3cm}
\begin{figure}[htbp]
\centering
\refstepcounter{algorithm}
\label{alg:bundle}

\begin{minipage}{0.99\columnwidth}
\hrule height 0.8pt
\vspace{0.25em}

\noindent\textbf{Algorithm \thealgorithm:} Three-cut Proximal Bundle Method

\vspace{0.15em}
\hrule height 0.4pt
\vspace{0.25em}

\noindent\textbf{Initialize} \(\mathbf{x}^0\) and the bundle described in
Sec.~\ref{subsec:bundle:procedure}.

\vspace{0.15em}
\noindent\textbf{Step \(k\):} \((k\ge0)\) Compute \(\mathbf{z}^{k+1}\) via
Algorithm~\ref{alg:prox_solver}.

\vspace{0.15em}
\noindent\hspace*{1em}\textbf{If} \eqref{eq:descent_step_test_letter} holds,
\(\mathbf{x}^{k+1}\leftarrow\mathbf{z}^{k+1}\);
\textbf{else} \(\mathbf{x}^{k+1}\leftarrow\mathbf{x}^{k}\).

\vspace{0.15em}
\noindent\hspace*{1em}\textbf{Stop:} if \(\Delta_k\le\varepsilon\), terminate.

\vspace{0.15em}
\noindent\hspace*{1.0em}Update \(\bar{\mathbf h}_k\), three-cut model, and
\(\underline f_{k+1}\) by \eqref{eq:agg_subgrad},
\eqref{eq:agg_cut_update_letter}, and \eqref{eq:three_cut_model}.

\vspace{0.25em}
\hrule height 0.8pt
\end{minipage}
\end{figure}

\vspace{-0.5cm}

\section{Numerical results}\label{sec:numerical}

\begin{figure}
\centering
\vspace{-1.35em}
\includegraphics[width=0.95 \columnwidth]{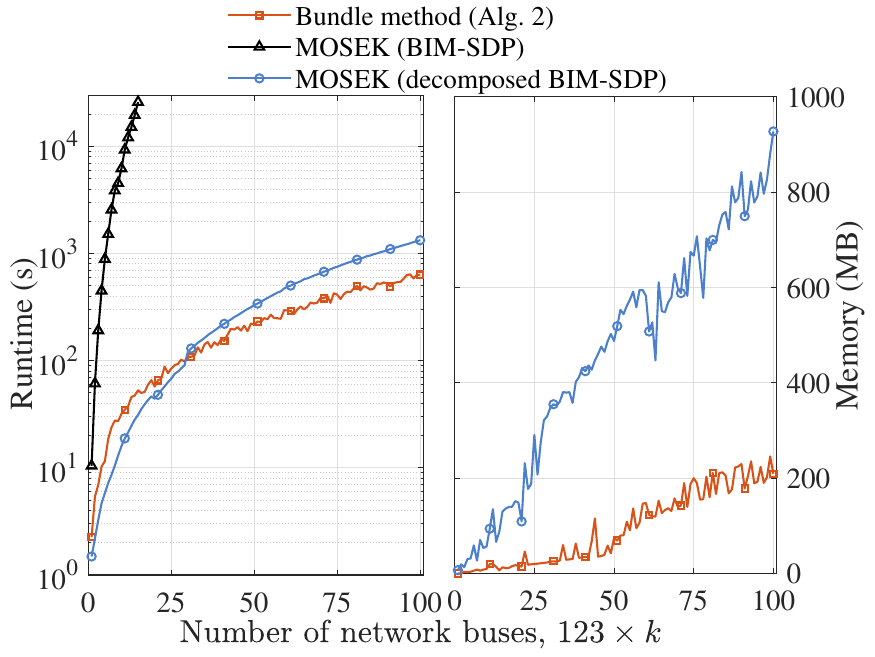}
\vspace{-1.35em}
\caption{Runtime and memory comparison versus network size with $k$ shown on the x-axis. The left panel shows runtime on a logarithmic scale for the bundle method, MOSEK, and MOSEK with the decomposed BIM-SDP, while the right panel shows memory usage for the bundle method and MOSEK with the decomposed BIM-SDP. The right panel excludes MOSEK with BIM-SDP because its memory usage becomes prohibitive as the network size grows; on the test platform, it fails at $k=16$ due to memory limitations. The simulations were conducted in MATLAB R2021b on a desktop computer with an Intel Core Ultra 7 265 processor at 2.40 GHz, 32 GB RAM, and MOSEK 11.1.}
\label{fig:runtime:memory}
\vspace{-0.5cm}
\end{figure}

We compare the bundle method
with MOSEK on the BIM-SDP in \eqref{eq:FP_compact_slack_SDP} and with MOSEK on a decomposed BIM-SDP
baseline adapted from \cite{StevenPartI}, which replaces the single large PSD constraint
by smaller per-line PSD constraints. 
We test on the IEEE 123$k$-bus test networks are built by connecting $k$ copies of the IEEE 123-bus feeder to a common root bus. 
As shown in Fig.~\ref{fig:runtime:memory}, MOSEK on the
BIM-SDP becomes memory prohibitive as the network size grows and is solvable
only up to $k=15$ in our tests; on this largest solvable case, the bundle
method is over $400\times$ faster. Compared with MOSEK on the decomposed
BIM-SDP, bundle method is about $2\times$ faster and uses about $75\%$ less memory at
the largest tested size, while matching the CVX/MOSEK benchmark within
$2\times 10^{-7}$ relative error.

For each iteration, Algorithm~\ref{alg:prox_solver} is at least
\(73\times\) faster than MOSEK across variety network sizes, with relative error
at most \(8.3\times10^{-5}\); see online Appendix~\ref{appen:sec:numerical} in details.

\vspace{-0.3cm}
\section{Conclusion}
This paper establishes an exact vectorized BIM-SDP reformulation for the three-phase unbalanced power flow and develops a scalable three-cut proximal
bundle solver for large-scale distribution network feasibility assessment.
Numerical results on IEEE test cases show improved runtime and memory
efficiency relative to MOSEK. Future work will extend the solver to
hosting-capacity and interconnection studies under diverse feeder topologies
and time-coupled constraints, while improving the current empirical $\rho$ with adaptive stepsize rules.


\appendices

\section{Dual of \texorpdfstring{\eqref{eq:FP_compact_slack_SDP}}{(1)}}\label{appen:dual:FPbeta}

Rather than directly penalizing the primal positive semidefinite (PSD) constraint in \eqref{eq:FP_compact_slack_SDP}, we derive the dual SDP of \eqref{eq:FP_compact_slack_SDP} for two reasons.
First, a primal penalty would retain the high-dimensional matrix variable $\mathbf{W}$, whereas the dual SDP is parameterized by the lower-dimensional variables $(\mathbf{y},\bm{\Gamma})$.
Second, at a rank-one optimal primal solution $\mathbf{W}^\star$, the matrix $-\mathbf{W}^\star$ has the largest eigenvalue $0$ with multiplicity $3N-1$. 
Hence, the subgradient associated with an analogous primal spectral penalty $\phi(\mathbf{W})=[\lambda_{\max}(-\mathbf{W})]_+$ would be highly nonunique (due to multiple eigenvectors), making the function $\phi$ highly nonsmooth near $\mathbf{W}^\star$. 
Such a highly nonsmooth structure near $\mathbf{W}^\star$ is difficult for the bundle method to capture. 
In contrast, the dual PSD matrix $\mathbf{H}^\star$ from \eqref{eq:dual_FP_betac} typically has only one zero eigenvalue. Thus, $-\mathbf{H}^\star$ has a simple largest eigenvalue at $0$, making the dual spectral penalty much less nonsmooth, which helps the bundle method capture the nonsmooth structure.
This motivates applying the exact penalty to the dual PSD constraint.

Introduce the dual variables \(\mathbf y\in\mathbb R^{18N}_+\) and
\(\bm\Gamma\in\mathbb H^3\) associated with
\eqref{eq:FP_compact_slack_SDPb} and
\eqref{eq:FP_compact_slack_SDPc}, respectively.
The Lagrangian of \eqref{eq:FP_compact_slack_SDP} is
\begin{equation}\label{eq:appen_Lagrangian}
\begin{aligned}
\mathcal L(\mathbf W,\mathbf z;\mathbf y,\bm\Gamma)
={}& \beta\cdot\mathbf 1^\top \mathbf z+\operatorname{tr}(\mathbf C\mathbf W) \\
&+\mathbf y^\top\!\big(\mathcal A(\mathbf W)+\mathbf m(\mathbf u)-\mathbf z\big) \\
&+\operatorname{tr}\!\big(\bm\Gamma(\mathcal B(\mathbf W)-\mathbf M_1)\big),
\end{aligned}
\end{equation}
where $\mathcal{B}:\mathbb{H}^{3N}\to\mathbb{H}^3$ denote the block-extraction operator $\mathcal{B}(\mathbf{W}) := [\mathbf{W}]_{11}$ and $\mathbf{M}_{1} := \mathbf{V}_{1}\mathbf{V}_{1}^{\mathsf H}$.

Define $\mathcal{A}^*$ and $\mathcal{B}^*$ as the adjoint operators of $\mathcal{A}$ and $\mathcal{B}$, respectively. In particular, $\mathcal{B}^*(\bm{\Gamma})$ places $\bm{\Gamma}$ in the top-left $3\times 3$ block of a $3N\times 3N$ Hermitian matrix, with zeros elsewhere. Using the adjoint identities
\[
\mathbf y^\top \mathcal A(\mathbf W)
=\operatorname{tr}\!\big(\mathcal A^*(\mathbf y)\mathbf W\big),
\quad
\operatorname{tr}\!\big(\bm\Gamma\,\mathcal B(\mathbf W)\big)
=\operatorname{tr}\!\big(\mathcal B^*(\bm\Gamma)\mathbf W\big),
\]
we rewrite the Lagrangian as
\begin{align*}
\mathcal L(\mathbf W,\mathbf z;\mathbf y,\bm\Gamma)
={}&
(\beta\cdot\mathbf 1-\mathbf y)^\top \mathbf z
+\mathbf m(\mathbf u)^\top \mathbf y
-\operatorname{tr}\!\big(\bm\Gamma\mathbf M_1\big)
\nonumber\\
&\quad
+\operatorname{tr}\!\Big(
\big[\mathbf C+\mathcal A^*(\mathbf y)+\mathcal B^*(\bm\Gamma)\big]\mathbf W
\Big).
\end{align*}

Hence, the dual function
\[
a(\mathbf y,\bm\Gamma)
=
\inf_{\mathbf W\succeq 0,\;\mathbf z\ge \mathbf 0}
\mathcal L(\mathbf W,\mathbf z;\mathbf y,\bm\Gamma)
\]
is finite if and only if
\[
\mathbf 0\le \mathbf y\le \beta\cdot\mathbf 1,
\qquad
\mathbf C+\mathcal A^*(\mathbf y)+\mathcal B^*(\bm\Gamma)\succeq 0;
\]
otherwise, \(q(\mathbf y,\bm\Gamma)=-\infty\).
Indeed, if either condition fails, one can let a component of \(\mathbf z\) or
\(\mathbf W=t\mathbf v\mathbf v^{\mathsf H}\) with \(t\to+\infty\) drive
\(\mathcal L\) to \(-\infty\).

Therefore, the dual function is
\[
a(\mathbf y,\bm\Gamma)
=
\mathbf m(\mathbf u)^\top\mathbf y
-\operatorname{tr}(\bm\Gamma\mathbf M_1),
\]
and the Lagrange dual problem is
\begin{equation}
\label{eq:dual_max_appendix}
\begin{aligned}
\max_{\mathbf y,\bm\Gamma}\quad
& a(\mathbf{y}, \bm{\Gamma}) \\
\text{s.t.}\quad
& \mathbf 0\le \mathbf y\le \beta\,\mathbf\cdot \mathbf{1}, \\
& \mathbf C+\mathcal A^*(\mathbf y)+\mathcal B^*(\bm\Gamma)\succeq 0.
\end{aligned}
\end{equation}
Negating the objective of \eqref{eq:dual_max_appendix} yields \eqref{eq:dual_FP_beta}.

\section{Proof of Proposition \ref{prop:dual_exact_penalty}}
\label{appen:proofproposition1}
The following proof adapts the argument of \cite{ding2023revisiting} from real symmetric matrices to complex Hermitian matrices under the notation of \eqref{eq:dual_FP_beta}--\eqref{eq:exact_penalty_obj}.
Based on the definition in \eqref{eq:set}, we first rewrite \eqref{eq:dual_FP_beta} as the following equivalent problem with a single scalar inequality:
\begin{subequations}\label{eq:dual_scalar}
\begin{align}
\min_{(\mathbf y,\bm\Gamma)\in\mathcal X}\quad
& q(\mathbf y,\bm\Gamma)
:= -\mathbf m(\mathbf u)^\top \mathbf y
+ \operatorname{tr}(\bm\Gamma \mathbf M_1),
\label{eq:dual_scalar_a}\\
\text{s.t.}\quad
& \psi(\mathbf y,\bm\Gamma)
:= \lambda_{\max}\!\big(-\mathbf H(\mathbf y,\bm\Gamma)\big)\le 0.
\label{eq:dual_scalar_b}
\end{align}
\end{subequations}
Indeed,
\[
\mathbf H(\mathbf y,\bm\Gamma)\succeq 0
\quad\Longleftrightarrow\quad
\lambda_{\max}\!\big(-\mathbf H(\mathbf y,\bm\Gamma)\big)\le 0.
\]

Since \(\phi(\mathbf y,\bm\Gamma)=[\psi(\mathbf y,\bm\Gamma)]_+\), where $[a]_+:=\max\{a,0\}$, problem \eqref{eq:exact_penalty_obj} can be written as
\[
\min_{(\mathbf y,\bm\Gamma)\in\mathcal X}
\Bigl\{
q(\mathbf y,\bm\Gamma)+\alpha\phi(\mathbf y,\bm\Gamma)
\Bigr\},
\]
where function $q$ is defined in equation \eqref{eq:dual_scalar_a}.

Let \(d^\star\) and \(d_\alpha^\star\) denote the optimal values of
\eqref{eq:dual_scalar} and \eqref{eq:exact_penalty_obj}, respectively.
Let \((\mathbf W^\star,\mathbf z^\star)\) be an optimal solution of
\eqref{eq:FP_compact_slack_SDP}.

\underline{Step 1: we first prove \(d_\alpha^\star \le d^\star\).}
For any feasible point \((\mathbf y,\bm\Gamma)\) of \eqref{eq:dual_scalar},
one has \(\psi(\mathbf y,\bm\Gamma)\le 0\), and hence
$\phi(\mathbf y,\bm\Gamma)
=0$.
Therefore, on the feasible
set of \eqref{eq:dual_scalar}, the objective of
\eqref{eq:exact_penalty_obj} reduces to the objective of
\eqref{eq:dual_scalar}:
\[
q(\mathbf y,\bm\Gamma)+\alpha\phi(\mathbf y,\bm\Gamma)
=
q(\mathbf y,\bm\Gamma).
\]
Since \eqref{eq:exact_penalty_obj} minimizes over the larger
set \(\mathcal X\), whereas \eqref{eq:dual_scalar} additionally
imposes \(\psi(\mathbf y,\bm\Gamma)\le 0\), we have
\[
d_\alpha^\star
\le
\min_{\substack{(\mathbf y,\bm\Gamma)\in\mathcal X\\
\psi(\mathbf y,\bm\Gamma)\le 0}}
\left\{
q(\mathbf y,\bm\Gamma)+\alpha\phi(\mathbf y,\bm\Gamma)
\right\}
=
d^\star .
\]

\underline{Step 2: we next prove \(d_\alpha^\star \ge d^\star\).}
Fix any \((\mathbf y,\bm\Gamma)\in\mathcal X\).
Since \(\mathbf W^\star\succeq 0\) and \(\operatorname{tr}(\mathbf W^\star)>0\), the matrix
\(\mathbf W^\star/\operatorname{tr}(\mathbf W^\star)\) is positive semidefinite and has unit trace. Hence, the Rayleigh--Ritz characterization gives
\begin{equation}\label{eq:appen_psi_lb}
\psi(\mathbf y,\bm\Gamma)
\ge
-\frac{\operatorname{tr}\!\big(\mathbf H(\mathbf y,\bm\Gamma)\mathbf W^\star\big)}
{\operatorname{tr}(\mathbf W^\star)}.
\end{equation}
Hence, given \(\alpha> \operatorname{tr}(\mathbf W^\star)\), by \eqref{eq:appen_psi_lb} we have
\begin{align}\nonumber
\alpha\phi(\mathbf y,\bm\Gamma)
\geq \operatorname{tr}(\mathbf W^\star)\psi(\mathbf y,\bm\Gamma)\geq
-\operatorname{tr}\!\big(\mathbf H(\mathbf y,\bm\Gamma)\mathbf W^\star\big),
\label{eq:appen_penalty_lb}
\end{align}
where the first inequality follows from $\phi=[\psi]_{+}\geq\psi$.

Therefore,
\begin{equation}\label{eq:appen_main_lb}
q(\mathbf y,\bm\Gamma)+\alpha\phi(\mathbf y,\bm\Gamma)
\ge
q(\mathbf y,\bm\Gamma)
-\operatorname{tr}\!\big(\mathbf H(\mathbf y,\bm\Gamma)\mathbf W^\star\big).
\end{equation}

We next derive a lower bound for the right-hand side of \eqref{eq:appen_main_lb}.
Using the definitions of \(q\) in equation \eqref{eq:dual_scalar_a} and \(\mathbf H\) in constraint \eqref{eq:dual_FP_betac}, we obtain
\begin{align}
&q(\mathbf y,\bm\Gamma)
-\operatorname{tr}\!\big(\mathbf H(\mathbf y,\bm\Gamma)\mathbf W^\star\big) \nonumber\\
&=
-\mathbf m(\mathbf u)^\top \mathbf y
+\operatorname{tr}(\bm\Gamma \mathbf M_1)
-\operatorname{tr}(\mathbf C\mathbf W^\star)
\nonumber\\
&\quad
-\operatorname{tr}\!\big(\mathcal A^*(\mathbf y)\mathbf W^\star\big)
-\operatorname{tr}\!\big(\mathcal B^*(\bm\Gamma)\mathbf W^\star\big)
\nonumber\\
&=
-\operatorname{tr}(\mathbf C\mathbf W^\star)
-\mathbf y^\top\!\bigl(\mathcal A(\mathbf W^\star)+\mathbf m(\mathbf u)\bigr)
\nonumber\\
&\quad
+\operatorname{tr}\!\bigl(\bm\Gamma(\mathbf M_1-\mathcal B(\mathbf W^\star))\bigr),
\label{eq:appen_expand1}
\end{align}
where the second equality uses the adjoint identities
\[
\operatorname{tr}\!\big(\mathcal A^*(\mathbf y)\mathbf W^\star\big)
\!=\!
\mathbf y^\top \mathcal A(\mathbf W^\star),
\operatorname{tr}\!\big(\mathcal B^*(\bm\Gamma)\mathbf W^\star\big)
\!=\!
\operatorname{tr}\!\big(\bm\Gamma\,\mathcal B(\mathbf W^\star)\big).
\]
Using
\(\mathcal B(\mathbf W^\star)=\mathbf M_1\), \eqref{eq:appen_expand1} reduces to
\begin{align}
q(\mathbf y,\bm\Gamma)
-\operatorname{tr}\!\big(\mathbf H(\mathbf y,\bm\Gamma)\mathbf W^\star\big)
&=
-\operatorname{tr}(\mathbf C\mathbf W^\star)
\nonumber\\
&\quad
-\mathbf y^\top\!\bigl(\mathcal A(\mathbf W^\star)+\mathbf m(\mathbf u)\bigr).
\label{eq:appen_expand2}
\end{align}
By primal feasibility of \eqref{eq:FP_compact_slack_SDP}, we have
\[
\mathcal A(\mathbf W^\star)+\mathbf m(\mathbf u)\le \mathbf z^\star,
\qquad
\mathbf z^\star\ge 0.
\]
Also, \((\mathbf y,\bm\Gamma)\in\mathcal X\) implies
\(0\le \mathbf y\le \beta\cdot\mathbf 1\). Hence
\[
-\mathbf y^\top\big(\mathcal A(\mathbf W^\star)+\mathbf m(\mathbf u)\big)
\ge
-\mathbf y^\top \mathbf z^\star
\ge
-\beta\cdot\mathbf 1^\top\mathbf z^\star .
\]
Substituting the above bound into \eqref{eq:appen_expand2} yields
\[
q(\mathbf y,\bm\Gamma)
-\operatorname{tr}\!\big(\mathbf H(\mathbf y,\bm\Gamma)\mathbf W^\star\big)
\ge
-\operatorname{tr}(\mathbf C\mathbf W^\star)
-\beta\cdot\mathbf 1^\top \mathbf z^\star=d^\star,
\]
where the equality follows from the strong duality for
\eqref{eq:FP_compact_slack_SDP}, since 
\eqref{eq:dual_scalar} is obtained from the Lagrange dual by changing the sign of the objective and writing it as a minimization problem.

Therefore, 
\begin{equation}\label{eq:dual:lowerbound}
    q(\mathbf y,\bm\Gamma)
-\operatorname{tr}\!\big(\mathbf H(\mathbf y,\bm\Gamma)\mathbf W^\star\big)
\ge d^\star.
\end{equation}

Combining \eqref{eq:dual:lowerbound} with \eqref{eq:appen_main_lb}, we obtain
\begin{equation}\label{eq:appen_main_lc}
q(\mathbf y,\bm\Gamma)+\alpha\phi(\mathbf y,\bm\Gamma)\ge d^\star,
\qquad
\forall (\mathbf y,\bm\Gamma)\in\mathcal X.
\end{equation}
Taking the minimum over $\mathcal X$ in \eqref{eq:appen_main_lc} yields $d_\alpha^\star \ge d^\star$.
Together with step 1, this gives $d_\alpha^\star=d^\star$.

\underline{Step 3: We prove the equality of optimal solution sets.}
Let \((\bar{\mathbf y},\bar{\bm\Gamma})\) be any minimizer of
\eqref{eq:exact_penalty_obj}. We show that it must satisfy
\(\psi(\bar{\mathbf y},\bar{\bm\Gamma})\le 0\).
Assume to the contrary that \(\psi(\bar{\mathbf y},\bar{\bm\Gamma})>0\).
Then \(\phi(\bar{\mathbf y},\bar{\bm\Gamma})=\psi(\bar{\mathbf y},\bar{\bm\Gamma})\), and if
\(\alpha>\operatorname{tr}(\mathbf W^\star)\), then by \eqref{eq:appen_psi_lb},
\begin{align*}
\alpha\phi(\bar{\mathbf y},\bar{\bm\Gamma})
\!=\!
\alpha\psi(\bar{\mathbf y},\bar{\bm\Gamma}) 
\!>\!
\operatorname{tr}(\mathbf W^\star)\psi(\bar{\mathbf y},\bar{\bm\Gamma}) 
\!\ge\!
-\operatorname{tr}\!\big(\mathbf H(\bar{\mathbf y},\bar{\bm\Gamma})\mathbf W^\star\big).
\end{align*}
Therefore,
\[
q(\bar{\mathbf y},\bar{\bm\Gamma})
+\alpha\phi(\bar{\mathbf y},\bar{\bm\Gamma})
>
q(\bar{\mathbf y},\bar{\bm\Gamma})
-\operatorname{tr}\!\big(\mathbf H(\bar{\mathbf y},\bar{\bm\Gamma})\mathbf W^\star\big)
\ge d^\star,
\]
where the second inequality follows from \eqref{eq:dual:lowerbound}. 
But this contradicts \(d_\alpha^\star=d^\star\). Hence every minimizer of
\eqref{eq:exact_penalty_obj} is feasible for \eqref{eq:dual_scalar}, and thus
optimal for \eqref{eq:dual_FP_beta}. 

The reverse inclusion is straightforward. 
Any minimizer of \eqref{eq:dual_scalar} is feasible for
\eqref{eq:exact_penalty_obj}, and the penalty term vanishes on the
feasible set of \eqref{eq:dual_scalar}. 
Since $d^\star=d_\alpha^\star$, minimizer of \eqref{eq:dual_scalar} also attains the optimal value of
\eqref{eq:exact_penalty_obj}. 
Thus, \eqref{eq:dual_scalar} and \eqref{eq:exact_penalty_obj}
have the same optimal solution set.

\section{Dual of the Proximal Subproblem and Recovery Formula}
\label{app:proof_prop_prox_dual}

This appendix derives the dual problem of \eqref{eq:prox_subproblem_epi} and the recovery
formula \eqref{eq:trial_point_general_short} stated in Proposition~\ref{prop:prox_dual}, where \(\boldsymbol\theta\) denotes the multiplier vector
associated with the three epigraph constraints in \eqref{eq:prox_subproblem_epi}. 
The optimality of \(\boldsymbol\theta\) returned by Algorithm~\ref{alg:prox_solver} is proved in Appendix~\ref{app:prox_cases}.

\subsection{Derivation of the Dual Problem}

The Lagrangian of \eqref{eq:prox_subproblem_epi} is
\begin{equation}\label{eq:app_lagrangian_prox}
\mathcal L(\mathbf x,r;\bm\theta)
=
r+\frac{\rho}{2}\|\mathbf x-\mathbf x^k\|^2
+\sum_{i=1}^3
\theta_i\bigl(a_k^i+\langle \mathbf h_k^i,\mathbf x\rangle-r\bigr),
\end{equation}
where $\bm\theta\in\mathbb R_+^3$ is the vector of dual multipliers associated with the
three epigraph constraints.

Collecting the terms involving \(r\), we have
\[
\mathcal L(\mathbf x,r;\bm\theta)
=
\Big(1-\sum_{i=1}^3\theta_i\Big)r
+
\frac{\rho}{2}\|\mathbf x-\mathbf x^k\|^2
+
\sum_{i=1}^3\theta_i
\big(a_k^i+\langle \mathbf h_k^i,\mathbf x\rangle\big).
\]
Since \(r\in\mathbb R\) is unconstrained, the infimum of $\mathcal{L}$ over \(r\) is finite
only if
$1-\sum_{i=1}^3\theta_i=0$.
Together with \(\bm\theta\in\mathbb R_+^3\), this gives
\(\bm\theta\in\Delta_3\).
Substituting this relation into \eqref{eq:app_lagrangian_prox} yields \eqref{eq:prox_dual_b},
and thus the dual problem \eqref{eq:prox_dual}.
 
\begin{subequations}\label{eq:prox_dual}
\begin{align}
\max_{\bm\theta\in\Delta_3}\quad & q_k(\bm\theta), \label{eq:prox_dual_a}\\
q_k(\bm\theta):={}&
\sum_{i=1}^3 \theta_i a_k^i
+\min_{\mathbf x\in\mathcal X}
\Biggl\{
\biggl\langle \sum_{i=1}^3 \theta_i \mathbf h_k^i,\mathbf x \biggr\rangle
+\frac{\rho}{2}\|\mathbf x-\mathbf x^k\|^2
\Biggr\}, \label{eq:prox_dual_b}
\end{align}
\end{subequations}
where $\Delta_3:=\{\bm\theta\in\mathbb R_+^3:\mathbf 1^\top\bm\theta=1\}$.

\subsection{Proof of the Recovery Formula}

Since \eqref{eq:prox_subproblem_epi} is a convex quadratic program and admits
a strictly feasible point, strong duality holds. For any fixed
$\bm\theta\in\Delta_3$, the inner minimization in \eqref{eq:prox_dual_b} is
strongly convex over $\mathcal X$, and hence has a unique minimizer
$\mathbf{z}^{k+1}(\bm\theta)$. 
First-order optimality condition for the inner minimization in \eqref{eq:prox_dual_b} is
\begin{equation}\label{eq:app_stationarity_x}
   \mathbf 0 \in
\rho\bigl(\mathbf z^{k+1}(\bm\theta)-\mathbf x^k\bigr)
+\sum_{i=1}^3 \theta_i \mathbf h_k^i
+\mathcal N_{\mathcal X}\bigl(\mathbf z^{k+1}(\bm\theta)\bigr), 
\end{equation}
where $\mathcal N_{\mathcal X}(\cdot)$ denotes the normal cone of $\mathcal X$.
By the characterization of Euclidean projection,
\eqref{eq:app_stationarity_x} is equivalent to
\[
\mathbf z^{k+1}(\bm\theta)
=
\Pi_{\mathcal X}\!\left(
\mathbf x^k-\frac{1}{\rho}\sum_{i=1}^3 \theta_i \mathbf h_k^i
\right).
\]
Applying this relation to an optimal dual solution
\(\boldsymbol\theta^k \in \arg\max_{\boldsymbol\theta\in\Delta_3} q_k(\boldsymbol\theta)\)
yields \eqref{eq:trial_point_general_short}.

Geometrically, $\Delta_3$ is a two-dimensional simplex parameterized by
$(\theta_1,\theta_2)$ with $\theta_3=1-\theta_1-\theta_2$.
A vertex has exactly one positive component, the relative interior of an
edge has exactly two positive components, and the interior of $\Delta_3$
has three positive components. 
These give the three cases considered in Appendix~\ref{app:prox_cases}.

\section{Sequential Case-by-Case Solver for the Proximal Subproblem}
\label{app:prox_cases}

This appendix presents the implementation details of the \emph{sequential} solver used in Algorithm~\ref{alg:prox_solver} and proves that Algorithm~\ref{alg:prox_solver} identifies an optimal solution
\(\boldsymbol\theta^\star\in\Delta_3\) of the dual problem \eqref{eq:prox_dual}.
Once such an optimal solution is obtained, the trial point \(\mathbf z^{k+1}\) follows from the
recovery formula \eqref{eq:trial_point_general_short} derived in Appendix~\ref{app:proof_prop_prox_dual}.

The case search uses the following KKT certificate. Since
\eqref{eq:prox_subproblem_epi} is a convex quadratic program and strong
duality holds, a candidate
\((\mathbf z^\star,r^\star,\boldsymbol\theta^\star)\) is optimal if and only if
\begin{subequations}\label{eq:app_kkt_certificate}
\begin{align}
&\mathbf z^\star
=\Pi_{\mathcal X}\left(
\mathbf x^k-\frac{1}{\rho}\sum_{i=1}^3
\theta_i^\star \mathbf h_k^i
\right),
\label{eq:app_kkt_projection}\\
&\boldsymbol\theta^\star\in\Delta_3,\qquad
\ell_k^i(\mathbf z^\star)\le r^\star,\quad i=1,2,3,
\label{eq:app_kkt_feasibility}\\
&\theta_i^\star
\bigl(\ell_k^i(\mathbf z^\star)-r^\star\bigr)=0,
\quad i=1,2,3 .
\label{eq:app_kkt_comp}
\end{align}
\end{subequations}
Because \(\boldsymbol\theta^\star\in\Delta_3\), at least one multiplier is
positive. Hence the KKT complementarity and feasibility imply
\[
r^\star=\max_{i=1,2,3}\ell_k^i(\mathbf z^\star),
\]
and any cut with a positive multiplier must be active at
\(\mathbf z^\star\).

Algorithm~\ref{alg:prox_solver} searches over the possible supports of
\(\boldsymbol\theta^\star\). It first tests the three vertices of
\(\Delta_3\), then the three edges, and finally the relative interior of
\(\Delta_3\). Each candidate is accepted once it satisfies
\eqref{eq:app_kkt_certificate}. This ordering avoids the two-dimensional
interior search whenever a lower-dimensional certificate is found.

\subsection{Case 1: Vertex of \texorpdfstring{$\Delta_3$}{Delta3}}

Suppose the multiplier \(\boldsymbol\theta^{k,\star}=\mathbf e_i\) is the vertex 
of \(\Delta_3\), where \(\mathbf e_i\) is the \(i\)-th unit vector in
\(\mathbb R^3\). By \eqref{eq:trial_point_general_short}, the corresponding
trial point is
\[
\mathbf z^{k+1}
=
\Pi_{\mathcal X}\left(
\mathbf x^k-\frac{1}{\rho}\mathbf h_k^i
\right).
\]

This candidate is accepted if the \(i\)-th cut attains the maximum model value
at \(\mathbf z^{k+1}\), namely,
\[
\ell_k^i(\mathbf z^{k+1})
\ge
\ell_k^j(\mathbf z^{k+1}),\qquad \forall j\ne i .
\]
If this test fails for all vertices, then the optimizer cannot lie at a
vertex, and the algorithm algorithm \emph{proceeds} to Case~2.

\subsection{Case 2: Edge of \texorpdfstring{$\Delta_3$}{Delta3}}\label{app:prox_cases:case2}

Consider the edge of \(\Delta_3\) spanned by \(\mathbf e_i\) and
\(\mathbf e_j\), with \(i\ne j\). A multiplier on this edge can be
parametrized as
\[
\theta_i=s,\
\theta_j=1-s,\
\theta_m=0,
\ m\in\{1,2,3\}\setminus\{i,j\},
\ s\in[0,1].
\]
The endpoints \(s=0\) and \(s=1\) correspond to the vertex cases already
checked in Case~1. Hence Case~2 only accepts edge candidates with
\(s\in(0,1)\). The associated recovered point is
\[
\mathbf z^{k+1}(s)
=
\Pi_{\mathcal X}\left(
\mathbf x^k-\frac{1}{\rho}
\left(s\mathbf h_k^i+(1-s)\mathbf h_k^j\right)
\right).
\]

For an edge candidate with support \(\{i,j\}\), i.e.,
\(\theta_i>0\) and \(\theta_j>0\), the complementarity condition
\eqref{eq:app_kkt_comp} implies that the two corresponding cuts must attain
the same value at the recovered point. Therefore, any accepted edge candidate
must satisfy
\begin{equation}
\varphi_{ij}(s)
:=
\ell_k^i\!\bigl(\mathbf z^{k+1}(s)\bigr)
-
\ell_k^j\!\bigl(\mathbf z^{k+1}(s)\bigr)
=0.
\label{eq:app_case2_equalheight}
\end{equation}

We next show that \(\varphi_{ij}\) is continuous, piecewise linear, and
\emph{nonincreasing} on \([0,1]\), which allows an exact left-to-right sweep over
the edge.

For each coordinate \(n\) in the box-constrained $\mathbf{y}$-block of \(\mathcal X\),
where \(x_n^k\), \(h_{k,n}^i\), and \(h_{k,n}^j\) denote the corresponding
entries of \(\mathbf x^k\), \(\mathbf h_k^i\), and \(\mathbf h_k^j\), the
pre-projection value is
\[
u_n(s)
=
x_n^k-\frac{1}{\rho}
\Bigl(
h_{k,n}^j+\bigl(h_{k,n}^i-h_{k,n}^j\bigr)s
\Bigr).
\]
A breakpoint occurs when \(u_n(s)\) hits one of the box bounds \(0\) or
\(\beta\). Whenever \(h_{k,n}^i-h_{k,n}^j\neq 0\), the corresponding
candidate breakpoints are
\[
s_{n,0}
=
\frac{\rho x_n^k-h_{k,n}^j}{\,h_{k,n}^i-h_{k,n}^j\,},
\qquad
s_{n,\beta}
=
\frac{\rho x_n^k-h_{k,n}^j-\rho\beta}{\,h_{k,n}^i-h_{k,n}^j\,}.
\]
Collecting all distinct values in \((0,1)\), together with \(0\) and \(1\),
yields a partition
\[
0=s_0<s_1<\cdots<s_P=1.
\]

On each open interval \((s_{p-1},s_p)\), the projection pattern is fixed.
Hence \(\mathbf z^{k+1}(s)\) is affine in \(s\) on that interval, and the
affine expression extends continuously to \([s_{p-1},s_p]\). Let
\(\mathcal F_p\) denote the set of coordinates on which the projection is
inactive on \((s_{p-1},s_p)\). For \(n\in\mathcal F_p\),
\[
z_n^{k+1}(s)=u_n(s),
\qquad
\frac{d}{ds}z_n^{k+1}(s)
=
-\frac{1}{\rho}\bigl(h_{k,n}^i-h_{k,n}^j\bigr),
\]
whereas the remaining coordinates are fixed at a bound and therefore have
zero derivative. Since
\[
\varphi_{ij}(s)
=
(a_k^i-a_k^j)
+
(\mathbf h_k^i-\mathbf h_k^j)^\top \mathbf z^{k+1}(s),
\]
we obtain, for \(s\in(s_{p-1},s_p)\),
\[
\begin{aligned}
\varphi_{ij}'(s)
&=
(\mathbf h_k^i-\mathbf h_k^j)^\top
\frac{d\mathbf z^{k+1}(s)}{ds} \\
&=
-\frac{1}{\rho}
\sum_{n\in\mathcal F_p}
\bigl(h_{k,n}^i-h_{k,n}^j\bigr)^2
=: -B_p
\le 0 .
\end{aligned}
\]

Thus \(\varphi_{ij}\) is nonincreasing on \([0,1]\). We now sweep the
intervals from left to right. The sweep is initialized by
\[
A_1:=\varphi_{ij}(s_0)=\varphi_{ij}(0).
\]
For \(p=1,\ldots,P\), suppose that the left-endpoint value
\(A_p=\varphi_{ij}(s_{p-1})\) is available. On
\([s_{p-1},s_p]\), \(\varphi_{ij}\) is given by
\[
\varphi_{ij}(s)
=
A_p-B_p(s-s_{p-1}),
\qquad s\in[s_{p-1},s_p].
\]
Define
\[
A_p^{\mathrm{end}}
:=
A_p-B_p(s_p-s_{p-1})
=
\varphi_{ij}(s_p).
\]

If \(A_p<0\), then, by monotonicity, no root can occur on this interval or on
any subsequent interval, and the current edge is discarded. If
\(A_p^{\mathrm{end}}>0\), then no root lies in \([s_{p-1},s_p]\), and the
sweep proceeds to the next interval with
\(A_{p+1}:=A_p^{\mathrm{end}}\).

It remains to consider the case \(A_p\ge 0\) and
\(A_p^{\mathrm{end}}\le 0\). If \(B_p>0\), then the unique root on this
interval is
\[
s^\star
=
s_{p-1}+\frac{A_p}{B_p}.
\]
If \(B_p=0\), then necessarily \(A_p=A_p^{\mathrm{end}}=0\), and
\eqref{eq:app_case2_equalheight} holds throughout the interval; in this
degenerate case, we select any
\(s^\star\in [s_{p-1},s_p]\cap(0,1)\), if such a point exists.

Given a selected root \(s^\star\in(0,1)\), the associated candidate
\(\mathbf z^{k+1}(s^\star)\) is accepted if the remaining cut does not
exceed the common active value, namely,
\begin{equation}
\ell_k^m\!\bigl(\mathbf z^{k+1}(s^\star)\bigr)
\le
\ell_k^i\!\bigl(\mathbf z^{k+1}(s^\star)\bigr)
=
\ell_k^j\!\bigl(\mathbf z^{k+1}(s^\star)\bigr),
\label{eq:app_case2_inactive}
\end{equation}
where \(m\) is the unique index in \(\{1,2,3\}\setminus\{i,j\}\).
If \eqref{eq:app_case2_inactive} fails, the current edge is rejected and
the next edge is tested. If no accepted edge candidate is found on the
three edges of \(\Delta_3\), the algorithm \emph{proceeds} to Case~3.

\subsection{Case 3: Interior of \texorpdfstring{$\Delta_3$}{Delta3}}\label{app:prox_cases:case3}
If no vertex or edge candidate is accepted, then any optimal
multiplier must lie in the relative interior of \(\Delta_3\). 
We use the
two-dimensional parametrization
\[
\mathcal T := \{(\theta_1,\theta_2): \theta_1\ge 0,\ \theta_2\ge 0,\ 
\theta_1+\theta_2\le 1\},
\]
with \(\theta_3=1-\theta_1-\theta_2\). We seek
\((\theta_1,\theta_2)\in \operatorname{int}(\mathcal T)\), so that
\(\theta_1,\theta_2,\theta_3>0\), and denote the corresponding recovered
point by \(\mathbf z^{k+1}(\theta_1,\theta_2)\). By \eqref{eq:app_kkt_comp},
all three cuts must attain the common epigraph value \(r\) at
\(\mathbf z^{k+1}(\theta_1,\theta_2)\). Therefore, the interior candidate
must satisfy
\begin{equation}
\mathbf F_k(\theta_1,\theta_2)
:=
\begin{bmatrix}
\ell_k^1\!\bigl(\mathbf z^{k+1}(\theta_1,\theta_2)\bigr)
-\ell_k^3\!\bigl(\mathbf z^{k+1}(\theta_1,\theta_2)\bigr)
\\[0.3em]
\ell_k^2\!\bigl(\mathbf z^{k+1}(\theta_1,\theta_2)\bigr)
-\ell_k^3\!\bigl(\mathbf z^{k+1}(\theta_1,\theta_2)\bigr)
\end{bmatrix}
=
\mathbf 0.
\label{eq:app_case3_F}
\end{equation}
Equivalently,
\begin{equation}\nonumber
\mathbf F_k(\theta_1,\theta_2)
=
\begin{bmatrix}
c_k^1+\langle \mathbf d_k^1,\mathbf z^{k+1}(\theta_1,\theta_2)\rangle\\
c_k^2+\langle \mathbf d_k^2,\mathbf z^{k+1}(\theta_1,\theta_2)\rangle
\end{bmatrix},
\label{eq:app_case3_F_cd}
\end{equation}
where
\[
c_k^i:=a_k^i-a_k^3,
\qquad
\mathbf d_k^i:=\mathbf h_k^i-\mathbf h_k^3,
\qquad i=1,2.
\]

By the projection formula \eqref{eq:trial_point_general_short},
\(\mathbf z^{k+1}(\boldsymbol\theta)\) is blockwise: the
\(\mathbf y\)-block is projected componentwise onto \([0,\beta]\), while the
\(\bm\Gamma\)-block is affine in \(\boldsymbol\theta\). Hence,
\(\mathbf F_k\) is nonsmooth only through the box projection on the
\(\mathbf y\)-block, and is piecewise affine and semismooth. We solve
\eqref{eq:app_case3_F} by a standard semismooth Newton method with
backtracking line search~\cite{Qi1993}.

At inner iteration $t$, let $(\theta_1^t,\theta_2^t)\in\operatorname{int}(\mathcal T)$
be the current iterate, and define the free index set on the $\mathbf y$-block by
\[
\mathcal F_t
:=
\{n:\ 0<[\mathbf y^{k+1}(\theta_1^t,\theta_2^t)]_n<\beta\}.
\]
Let $\mathbf D_t=\operatorname{diag}(m_n^t)$, where
\[
m_n^t=
\begin{cases}
1, & n\in\mathcal F_t,\\
0, & n\notin\mathcal F_t.
\end{cases}
\]

Accordingly, a generalized derivative of the blockwise projection is
represented by
\[
\mathbf P_t :=
\begin{bmatrix}
\mathbf D_t & 0\\
0 & \mathbf I_3
\end{bmatrix},
\]
where \(\mathbf I_3\) acts on the \(3\times 3\) \(\bm\Gamma\)-block.
Then a valid generalized Jacobian element of \(\mathbf F_k\) at
\((\theta_1^t,\theta_2^t)\) is
\[
\mathbf J_t
=
-\frac{1}{\rho}
\begin{bmatrix}
(\mathbf d_k^1)^\top \mathbf P_t \mathbf d_k^1
&
(\mathbf d_k^1)^\top \mathbf P_t \mathbf d_k^2\\
(\mathbf d_k^2)^\top \mathbf P_t \mathbf d_k^1
&
(\mathbf d_k^2)^\top \mathbf P_t \mathbf d_k^2
\end{bmatrix}.
\]

The semismooth Newton direction \(\mathbf p_t\) is obtained from
\[
\mathbf J_t \mathbf p_t = -\mathbf F_k(\theta_1^t,\theta_2^t).
\]
We then update
\[
(\theta_1^{t+1},\theta_2^{t+1})
=
(\theta_1^t,\theta_2^t)+\tau_t\mathbf p_t,
\]
where \(\tau_t\in(0,1]\) is chosen by backtracking line search so that
\((\theta_1^{t+1},\theta_2^{t+1})\in\operatorname{int}(\mathcal T)\)
and \(\|\mathbf F_k(\theta_1^{t+1},\theta_2^{t+1})\|_2\) is sufficiently
reduced. Since the current iterate lies in \(\operatorname{int}(\mathcal T)\),
sufficiently small \(\tau_t\) preserves the interior feasibility. Upon
convergence, the resulting \((\theta_1^\star,\theta_2^\star)\) yields the
desired trial point \(\mathbf z^{k+1}(\theta_1^\star,\theta_2^\star)\).

\section{Additional Simulation Results, Implementation Details, and Convergence Remarks}
\label{appen:sec:numerical}

Table~\ref{tab:subprob_time} compares the average time for computing the trial point $\mathbf z^{k+1}$ by
Algorithm~\ref{alg:prox_solver} and by MOSEK on the IEEE 123$k$-bus test networks,
for $k=1,10,50,100$. These four cases are selected to represent different
network scales. The results show that the per-iteration cost of
Algorithm~\ref{alg:prox_solver} remains substantially lower across all four
cases, and the runtime gap widens as the network size increases.

\begin{table}[t]
\caption{Average time (s) for computing $\mathbf z^{k+1}$ at four network scales}
\label{tab:subprob_time}
\centering
\renewcommand{\arraystretch}{1.05}
\begin{tabular}{ccccc}
\toprule
$k$ & 1 & 10 & 50 & 100 \\
\midrule
Alg.~\ref{alg:prox_solver} & $1.2\times10^{-3}$ & $3.7\times10^{-3}$ & $1.6\times10^{-2}$ & $5.7\times10^{-2}$ \\
MOSEK & $8.7\times10^{-2}$ & $3.1\times10^{-1}$ & $2.1$ & $4.7$ \\
\bottomrule
\end{tabular}
\end{table}

We set the initialization $\mathbf{x}^0:=\bigl(\tfrac{\beta}{2}\mathbf 1,\mathbf 0_3\bigr)\in\mathcal X$.
The slack-penalty parameter $\beta$ in \eqref{eq:FP_compact_slack_SDPa} is set to $0.1$, following
the numerical setting in \cite{fang2025}, under which the SDP relaxation was observed
to be exact for the tested radial networks. 
The penalty parameter in \eqref{eq:exact_penalty_obj}
is chosen as $\alpha = 2\cdot\mathbf 1^\top \overline{\mathbf v}$, which satisfies the condition $\alpha>\operatorname{tr}(\mathbf{W}^\star)$
in Proposition~\ref{prop:dual_exact_penalty}.
The proximal parameter $\rho$ in \eqref{eq:prox_subproblem_epi}, which is also used in Algorithm~\ref{alg:prox_solver}, is set to $4$. The serious-step parameter $\eta$ in \eqref{eq:descent_step_test_letter} is set to $0.1$, and the stopping tolerance $\varepsilon$ in Algorithm~\ref{alg:bundle} is set to $10^{-5}$.

By the convergence result in \cite{Diaz2023} and the exact-penalty equivalence in Proposition~\ref{prop:dual_exact_penalty}, Algorithm~\ref{alg:bundle} converges to an optimal solution of the dual problem \eqref{eq:dual_FP_beta}.


\begin{thebibliography}{1}
\providecommand{\url}[1]{#1}
\csname url@samestyle\endcsname
\providecommand{\newblock}{\relax}
\providecommand{\bibinfo}[2]{#2}
\providecommand{\BIBentrySTDinterwordspacing}{\spaceskip=0pt\relax}
\providecommand{\BIBentryALTinterwordstretchfactor}{4}
\providecommand{\BIBentryALTinterwordspacing}{\spaceskip=\fontdimen2\font plus
\BIBentryALTinterwordstretchfactor\fontdimen3\font minus \fontdimen4\font\relax}
\providecommand{\BIBforeignlanguage}[2]{{%
\expandafter\ifx\csname l@#1\endcsname\relax
\typeout{** WARNING: IEEEtran.bst: No hyphenation pattern has been}%
\typeout{** loaded for the language `#1'. Using the pattern for}%
\typeout{** the default language instead.}%
\else
\language=\csname l@#1\endcsname
\fi
#2}}
\providecommand{\BIBdecl}{\relax}
\BIBdecl

\bibitem{cui2019solvability}
B.~Cui and X.~A. Sun, ``Solvability of power flow equations through existence and uniqueness of complex fixed point,'' \emph{arXiv preprint arXiv:1904.08855}, 2019.

\bibitem{Molzahn2026DOE}
H.~Moring, D.~K. Molzahn, and J.~L. Mathieu, ``Exploring the inexactness of second-order cone relaxations for calculating operating envelopes,'' \emph{IEEE Trans. Power Syst.}, pp. 1--12, 2026.

\bibitem{fang2025}
B.~Fang, Y.~Chen, and C.~Zhao, ``Approximating dispatchable regions in three-phase radial networks with conditions for exact {SDP} relaxation,'' \emph{IEEE Trans. Power Syst.}, 2026, to be published; arXiv:2503.22385.

\bibitem{StevenPartI}
S.~H. Low, ``Convex relaxation of optimal power flow---{Part} {I}: Formulations and equivalence,'' \emph{IEEE Trans. Control Netw. Syst.}, vol.~1, no.~1, pp. 15--27, 2014.

\bibitem{ding2023revisiting}
L.~Ding and B.~Grimmer, ``Revisiting spectral bundle methods: Primal-dual (sub)linear convergence rates,'' \emph{SIAM J. Optim.}, vol.~33, no.~2, pp. 1305--1332, 2023.

\bibitem{Qi1993}
L.~Qi and J.~Sun, ``A nonsmooth version of {Newton's} method,'' \emph{Math. Program.}, vol.~58, no.~1, pp. 353--367, 1993.

\bibitem{Diaz2023}
M.~D{\'i}az and B.~Grimmer, ``Optimal convergence rates for the proximal bundle method,'' \emph{SIAM J. Optim.}, vol.~33, no.~2, pp. 424--454, 2023.

\end{thebibliography}
\end{document}